\documentstyle[psfig,12pt]{amsart}

\newcommand{\C}{\mathbb C}

\newcommand{\Z}{\mathbb Z}
\renewcommand{\pf}{{\bf Proof. }}
\newcommand{\lam}{\lambda}

\newcommand{\sig}{\sigma}

\renewcommand{\epsilon}{\varepsilon}
\renewcommand{\phi}{\varphi}
\newtheorem{lem}{Lemma}
\newtheorem{T}[lem]{Theorem}

\newtheorem{defi}{Definition}
\theoremstyle{definition}

\begin{document}

\title[ Braid Group Representations    ]{ Dimension $n$ Representations 
of the Braid Group on $n$ Strings      }
\author[I. Sysoeva        ]{Inna Sysoeva     }
\date{\today}
\address{Department of Mathematics \\
The Pennsylvania State University \\
University Park, PA 16802 }    
\email{sysoeva@@math.psu.edu   }



\maketitle

\section{Introduction}
Let  $B_n$ be the braid group
on $n$ strings. In his paper \cite{form} Ed Formanek classified all irreducible representations of  $B_n$ of dimension at most $(n-1).$ Since then there were some attempts to classify irreducible representations of $B_n$ of dimension $n$.
In particular the classification is known for very small $n$. Case $n=3$ was 
done by Ed Formanek (\cite{form}, Theorem 24). Woo Lee
has classified the four-dimensional irreducible representations of $B_4$ (\cite{WooLee}).

In this paper we solve this problem completely for $n\ge 9$. Before stating our main classification theorem let us describe the following representation of $B_n$ of dimension $n$.

\begin{defi} \label{defi:1}{\bf The standard representation} is the representation
$$\tau_n:B_n \to GL_n ({\Z}[t ^{\pm 1}]$$ defined by

\vskip 1cm
$$\rho (\sig_i)=\left( \begin{array}{ccccc}
I_{i-1}&&&\\
&0&t&\\
&1&0&\\
&&&I_{n-1-i}
\end{array} \right),$$ 
\vskip .5cm
\noindent
for $i=1,2,\dots, n-1,$ where $I_{k}$ is the $k\times k$ identity matrix.

\end{defi}

We call the above representation standard because of its simplicity. Surprisingly, it does not seem to be well-known. In fact it looks like it 
was first discovered only in 1996 by Dian-Min Tong, Shan-De Yang and Zhong-Qi Ma ( \cite{TongYangMa}, Equation (19)).

\begin{T}\label{T:final} Suppose that $\rho:B_{n} \to GL_{n}({\C})$ is an irreducible representation of $B_{n}$ of dimension $n\ge 9.$ Then it is a tensor product of a one-dimensional representation and a specialization for $u\neq 0, 1$ of the standard representation.

\end{T}

To explain the ideas of the proof, we need the following definition.

\begin{defi}  Suppose $\rho $ is a representation of the Artin braid group $B_n. $ The corank of $\rho $ is the rank of $\rho (\sigma_i),$ where $\sigma_i$ are the standard generators of $B_n.$ (This makes sense because all $\sigma_i$ are conjugate.)

\end{defi}

If one looks at the proof of the classification theorem of Formanek in \cite{form}, it can be separated into two parts. The first is to classify all irreducible representations of braid groups of corank $1.$ The second is to prove that apart from a few exceptions, the irreducible  representations of braid groups $B_n$ of dimension at most $(n-1)$ can be obtained as a tensor product of a one-dimensional representation and an irreducible representation of corank $1.$

Our proof follows a similar strategy. The first part of it, the classification of irreducible representations of corank $2$ was carried out in \cite{tes}. In this paper we complete the proof of Theorem \ref{T:final} by proving that for $n\ge 9$ every irreducible representation of $B_{n}$ of dimension $n$ is the tensor product of a one-dimensional representation and a representation of corank
$2$.

{\bf Acknowledgments:} The author wishes to express her deep gratitude to Ed Formanek for numerous helpful discussions and support. His remarks simplified considerably the proof of the main theorem. In particular, Lemma \ref{lem:size} is due to him.

\section {Proof of Theorem \ref{T:final}}

We proved in \cite{tes}, Theorem 5.5 and Corollary 5.6 that for $n\ge 7$ every irreducible complex representation of $B_n$ of corank $2$ is a specialization of the standard representation (see Definition \ref{defi:1}.) So to complete the proof of Theorem \ref{T:final} it is enough to show that  for $n\ge 9$ every irreducible representation of $B_{n}$ of dimension $n$ is the tensor product of a one-dimensional representation of corank $2$. This will be done in Theorem \ref{T:84}. Before that we need some preparatory results. The key of the proof is the following theorem, which is similar to Theorem 16 of \cite{form}.

\begin{T}\label{T:osn}
Suppose that $\rho:B_{n+1} \to GL_{n+1}({\C})$ is an irreducible
 representation of $B_{n+1}$ of dimension $n+1$
($n \geq 4$). Suppose that the restriction of $\rho,$
$\rho | {B_{n-1} \times <\sig_{n}  >},$ stabilizes one-dimensional subspace
${\C}v$ of $
{\C}^{n+1}.$ 

Then $rank(\rho(\sig_1)-yI)=2$ for some $y \in {\C}^{*}.$ 

\end{T}
\pf For notational simplicity we will write $\sig$ instead of
$\rho(\sig)$ for $\sig \in B_n.$

By hypothesis,
$$\rho | {B_{n-1} \times <\sig_n  >}:{\C}v \to {\C}v$$ is a one-dimensional
representation of $B_{n-1}\times B_2,$ so there exist $x,y \in {\C}^{*}$
such that
$$\sig_1 v=\sig_2 v= \dots =\sig_{n-2}v=yv, \,\,\, \sig_n v=xv$$

Consider $ \theta=\theta _{n+1}=\sig_1 \sig_2 \dots \sig_n,$
$\sig_0= \theta \sig_n \theta ^{-1},$
$$v_n=v,\,\, v_{n+1}=\theta v, \,\, 
v_1=\theta ^2 v, \dots, v_{n-1}=\theta ^n v.$$ Conjugation by 
$\theta$ permutes $\sig_1, \dots , \sig_n, \sig_0$ cyclically.

Because $\rho$ is an irreducible representation and $\theta ^{n+1}$
is central in $B_{n+1},$ $\rho(\theta ^{n+1})=dI$ for some $d \in {\C}^{*}.$
Thus, the left action of $\theta$ permutes ${\C}v_1,{\C}v_2, \dots , {\C}v_{n+1}$ cyclically.

We have:$$\sig_iv_i=xv_i,$$
$$\sig_iv_{i+j}=yv_{i+j}$$ for $$i=1,\dots ,n+1,\,\,\, j=2, \dots ,n-1,$$ where indices are taken modulo $n+1.$

The following table summarizes the above calculations:
\vskip .5cm

\begin{tabular}{c|c|c|c|c|c|c|c|c|c|c}

 &  & $v_1$ & $v_2$ & $v_3$ &$\dots$ & $v_{n-1}$ & $v_n$ & $v_{n+1}$  \\
\hline
 
 &  $\sig_1 $  & $xv_1$ &  & $yv_3$ & $\dots$& $yv_{n-1}$ & $yv_{n}$ &    \\
\hline
  &  $\sig_2 $  &  & $xv_2$ & & $\dots$& $yv_{n-1}$ & $yv_{n}$ & $yv_{n+1}$    \\
\hline
  &  $\sig_3 $  & $yv_1$ & &$xv_3$  & $\dots$& $yv_{n-1}$ & $yv_{n}$ & $yv_{n+1}$    \\
\hline

 & $\vdots$  &$\vdots$ &$\vdots$  &$\vdots$ & $\ddots$ & $\vdots$    &$\vdots$    & $\vdots$    \\
\hline
 
  &  $\sig_{n-1} $  & $yv_1$ & $yv_2$  &$yv_3$  & $\dots$& $xv_{n-1}$ &   & $yv_{n+1}$    \\
\hline
 &  $\sig_n $  & $yv_1$ & $yv_2$  & $yv_3$  & $\dots$& & $xv_{n}$ & \\
\hline 
&  $\sig_0 $  & $yv_1$ & $yv_2$  &$yv_3$  & $\dots$& $yv_{n-1}$ & & $xv_{n+1}$    \\

\end{tabular}

\vskip .5cm

Suppose that $v_1, \dots , v_{n+1}$ are linearly dependent.
Consider
$$a_1v_1+a_2v_2+ \dots + a_tv_t=a_1v_1+a_2 \theta v_1+ \dots +a_t{\theta}^{t-1}v_1=0,$$
a linear dependence relationship with minimal $t.$

In the equation above, $a_1 \neq 0,$ since $\theta$ is invertible, and $a_t \neq 0$ by the
minimality of $t.$

We claim that $t \geq n.$ Indeed, suppose that $t \leq n-1.$
Then $v_{n-1}$ is a linear combination of $v_1, \dots , v_{n-2},$
which are eigenvectors for $\sig_n$ with $\sig_nv_i=yv_i,\,\,\, i=1, \dots ,n-2.$
So, $\sig_n v_{n-1}=yv_{n-1}.$ Applying ${\theta}^3$ implies that
$\sig_2 v_1=yv_1,$ which means that ${\C}v_1$ is $B_{n+1}-$invariant,
which contradicts the irreducibility of $\rho.$ So, $t \geq n.$

Thus, $v_1, \dots , v_{n-1}$ are linearly independent.

Assume that $rank(\sig_1-yI)> 2.$ Then, as 
$$dim(Ker(\sig_1-yI))+ rank(\sig_1-yI)=n+1,$$
$dim(Ker(\sig_1-yI))\leq n-2.$ 

Note that $v_3, \dots ,v_n$ are $n-2$ linearly independent elements
of $L=Ker(\sig_1-yI).$ So, $dim(Ker(\sig_1-yI))= n-2,$
and $L=span \{v_3,\dots, v_n\}.$

Since vectors $\{v_1, \dots ,v_{n-1} \}$ are linearly independent,
$\{v_2, \dots ,v_{n} \}$ and $\{v_3, \dots ,v_{n+1} \}$ are
also linearly independent. Therefore $v_2 \notin L,$ and 
$v_{n+1} \notin L.$

The action of $\theta$ implies that for $i=1, \dots , n+1$
$$Ker(\sig_i-yI)=span \{ v_{i+2}, \dots , v_{i-2}\},$$
 $v_{i-1} \notin L,$ and $v_{i+1} \notin L,$ where 
indices are taken modulo $n+1.$

$\sig_1$ commutes with $\sig_n,$ and $n \geq 4,$
so
$$(\sig_n -yI)\sig_1v_2=\sig_1(\sig_n -yI)v_2=0.$$

Thus, $\sig_1v_2\in Ker (\sig_n -yI),$ so
$$\sig_1v_2=b_1v_1+b_2v_2 +
\dots + b_s v_s,$$ where $1 \leq s \leq n-2$ and $b_s \neq 0.$

We claim that $s\leq 2.$ Indeed, if $s \geq 3,$ then
$$0=\sig_1(\sig_{s+1}-yI)v_2=(\sig_{s+1}-yI)\sig_1v_2=$$
$$=(\sig_{s+1}-yI)(b_1v_1+b_2v_2 +
\dots + b_s v_s)=(\sig_{s+1}-yI)b_sv_s.$$ 
This contradicts  the fact that $v_{s} \notin Ker(\sig_{s+1}-yI).$

Thus, $$\sig_1v_2=b_1v_1+b_2v_2 , \,\,\, b_1,b_2 \in {\C}.$$

By a symmetric argument which reverses the roles of $\sig_1$ and $\sig_n,$
and starts with the equation
$$(\sig_1 -yI)\sig_nv_{n-1}=\sig_n(\sig_1 -yI)v_{n-1}=0,$$
we obtain 
$$\sig_nv_{n-1}=c_1v_{n-1}+c_2v_n , \,\,\, c_1,c_2 \in {\C}.$$

Using the action of $\theta,$ we get the following table:
\vskip .5cm
\begin{tabular}{c|c|c|c|c|c|c|c|c|c|c}

 &  & $v_1$ & $v_2$ & $\dots$ & $v_n$ & $v_{n+1}$  \\
\hline
 
 &  $\sig_1 $  & $xv_1$ & $b_1v_1+b_2v_2$& $\dots$& $yv_{n}$ & $c_1v_{n+1}+c_2v_1 $  \\
\hline
  &  $\sig_2 $  & $c_1v_{ 1}+c_2v_2 $ & $xv_2$ & $\dots$& $yv_{n}$ & $yv_{n+1}$    \\
\hline
  &  $\sig_3 $  & $yv_1$ &$c_1v_{ 2}+c_2v_3 $  & $\dots$& $yv_{n}$ & $yv_{n+1}$    \\
\hline

 & $\vdots$  &$\vdots$ &$\vdots$  & $\ddots$ & $\vdots$    & $\vdots$    \\
\hline
 
  &  $\sig_{n-1} $  & $yv_1$ & $yv_2$  & $\dots$& $b_1v_{n-1}+b_2v_n $   & $yv_{n+1}$    \\
\hline
 &  $\sig_n $  & $yv_1$ & $yv_2$  & $\dots$&  $xv_{n}$ &$b_1v_{n}+b_2v_{n+1} $    \\

\end{tabular}
\vskip .5cm

$Span\{v_1, v_2, \dots , v_{n+1}\}$ is $B_{n+1}-$invariant.
Thus, if $\{v_1, v_2, \dots , v_{n+1}\}$ are linearly 
dependent, then $\rho$ is reducible. So,$\{v_1, v_2, \dots , v_{n+1}\}$ are linearly 
independent,  and they form a basis for ${\C}^{n+1}.$

In this basis:
$$\sig_1=\left( \begin{array}{cccccc}
x&b_1&&&&c_2\\
0&b_2&&&&\\
&&y&&&\\
&&&\ddots&&\\
&&&&y&\\
&&&&&c_1
\end{array} \right),\sig_3=\left( \begin{array}{ccccccc}
y&&&&&&\\
&c_1&&&&&\\
&c_2&x&b_1&&&\\
&&&b_2&&&\\
&&&&y&&\\
&&&&&\ddots&\\
&&&&&&y
\end{array} \right).$$

Using the $(3,2)-$entry of the matrix $\sig_1\sig_3=\sig_3\sig_1,$
we have
$$b_2c_2=yc_2.$$

If $c_2=0,$ then ${\C}v_1$ is invariant under $B_{n+1},$
which contradicts the irreducibility of $\rho.$ So, $c_2 \neq 0.$ 
Thus, $b_2=y.$ Then $rank (\sig_1-yI)\leq 2,$ a contradiction.

So, $rank (\sig_1-yI)\leq 2.$  But by \cite{form}, Theorem 10,
the case $rank (\sig_1-yI)=1$  is impossible. Thus, $rank (\sig_1-yI)=2.$ 
 \vskip .5cm

The following argument is due to E. Formanek. He also used it in \cite{form},
Lemma 17 and Corollary 18. My original argument was much longer. 

The next Lemma
\ref{lem:82} is a corollary of Theorem 23 of \cite{form}, which
classifies the irreducible representations of $B_n$ of dimension at most $n-1.$

\begin{lem}\label{lem:82} If $\rho:B_n \to GL_r ({\C})$ is irreducible and
$r \leq n-3,$ then $\rho$ is one-dimensional.

\end{lem}

\begin{lem}\label{lem:size} Let $\rho :B_n \to GL_r ({\C})$ be a representation,
where $n \geq 6.$ Suppose that $\lam$ is an eigenvalue of $\rho(\sig_{n-1}).$
Suppose that the largest Jordan block corresponding to $\lam$ has size $s$ and
 multiplicity $d.$ 

If $d \leq n-5,$ then $\rho|{B_{n-2} \times <\sig_{n-1}>}$ has
a one-dimensional invariant subspace.

\end{lem}
\pf Let $f(t)$ be the minimal polynomial of $\rho(\sig_{n-1}).$
Set $m(t)=f(t)/(t-\lam).$ Let $V$ be the image of ${\C}^r$ under
$m(\rho(\sig_{n-1})).$ Then $V$ is invariant under $\rho|{B_{n-2} \times <\sig_{n-1}>},$ and $dim V=d.$ If $d \leq n-5,$ then by Lemma \ref{lem:82},
all composition factors of $$\rho|{B_{n-2} \times <\sig_{n-1}>}:V \to V$$ 
are one-dimensional.

\begin{T} \label{T:84} For $n\geq 9,$ every $n-$dimensional complex irreducible representation $\rho$ of
the braid group $B_n$ is equivalent to
a tensor product of a one-dimensional representation 
$\chi (y),$ $y \in {\C}^{*},$ and an $n$-dimensional
representation of corank $2.$
\end{T}
\pf Assume not. Then by  Theorem \ref{T:osn} and Lemma \ref{lem:size},
the largest Jordan block corresponding to every eigenvalue of $\rho(\sig_{n-1})$
has multiplicity $\geq n-4.$

If $\rho(\sig_{n-1})$ has two or more eigenvalues, we get $$(n-4)+(n-4) \leq n,$$
a contradiction, since $n \geq 9.$ Similarly, if some eigenvalue
has the corresponding largest Jordan block of size $s \geq 2,$ we get a contradiction $$2(n-4) \leq n.$$

Thus, $\rho(\sig_{n-1})$ has only one eigenvalue $\lam$ and the Jordan canonical
form of $\rho(\sig_{n-1})$ consists of $1 \times 1$ elementary
Jordan blocks. But then $\rho(\sig_{n-1})=\lam I,$ 
which contradicts the irreducibility of $\rho.$

This completes the proof of the theorem, and thus the proof of Theorem \ref{T:final}.

\end{document}